\documentclass[12pt,reqno]{amsart}
\usepackage{amssymb}
\usepackage{amsmath}
\usepackage{mathrsfs}

\usepackage{amssymb, amsmath, amsfonts, latexsym}
\usepackage{enumerate}

\newtheorem{thm}{Theorem}[section]
\newtheorem{cor}[thm]{Corollary}
\newtheorem{lem}[thm]{Lemma}
\newtheorem{prop}[thm]{Proposition}

\newtheorem{remarks}[thm]{Remark}

\theoremstyle{definition}
\newtheorem{defn}{Definition}[section]

\numberwithin{equation}{section} \theoremstyle{remark}

\title[CLT and Moderate deviations for a  class of semilinear SPDEs]{Central limit theorem and moderate deviations for a  class of semilinear SPDEs}

\author{Shulan Hu}
\address{Shulan Hu \\School of Statistics and Mathematics, Zhongnan University of Economics and Law, 430073, PR China.}
\email{hu\_shulan@yahoo.com}

\author{Ruinan Li}
\address{Ruinan Li \\School of Statistics and Information, Shanghai University of International Business and Economics, Shanghai, 201620, PR China.}
\email{ruinanli@amss.ac.cn}

\author{Xinyu Wang}
\address{Xinyu Wang\\ School of Mathematics and Statistics, Huazhong University of Science and Technology, 430074, PR China.}
\email{wang\_xin\_yu@hust.edu.cn}

\newcommand{\ee}{\mathbb{E}}

\newcommand{\rr}{\mathbb{R}}
\newcommand{\pp}{\mathbb{P}}

\def\AA{\mathcal A}

\def\FF{\mathcal F}
\def\EE{\mathcal E}

\def\HH{\mathcal H}

\def\e{\varepsilon}

\def\<{\langle}
\def\>{\rangle}

\def\beq{\begin{equation}}
\def\nneq{\end{equation}}

\def\bdef{\begin{defn}}
\def\ndef{\end{defn}}

\def\bthm{\begin{thm}}
\def\nthm{\end{thm}}

\def\bprop{\begin{prop}}
\def\nprop{\end{prop}}

\def\brmk{\begin{remarks}}
\def\nrmk{\end{remarks}}

\def\bexa{\begin{exa}}
\def\nexa{\end{exa}}

\def\blem{\begin{lem}}
\def\nlem{\end{lem}}

\def\bcor{\begin{cor}}
\def\ncor{\end{cor}}

\date{}

\def\bexe{\begin{exe}}
\def\nexe{\end{exe}}

\def\bprf{\begin{proof}}
\def\nprf{\end{proof}}

\def\bdes{\begin{description}}
\def\ndes{\end{description}}

\begin{document}
\maketitle

\begin{abstract} In this paper we prove  a central limit theorem and a  moderate deviation principle for  a class of semilinear stochastic partial differential equations, which  contain Burgers' equation and the stochastic reaction-diffusion equation. The weak convergence method plays an important role.
\end{abstract}
\vskip 0.4cm
\noindent {\bf MSC 2010: } Primary 60H15; Secondary 60F05, 60F10.
\vskip 0.4cm

\noindent {\bf keywords:} Burgers' equation, Stochastic reaction-diffusion equation, Large deviations,  Moderate deviations.

\section{Introduction}

 For any $\e>0$,  consider the following semilinear stochastic partial differential equation
\begin{align}\label{SPDE}
       \frac{\partial U^{\e}}{\partial t}(t,x)=&\frac{\partial^2 U^{\e}}{\partial x^2}(t,x)+\sqrt{\e}\sigma(t,x,U^{\e}(t,x))\frac{\partial ^2W}{\partial t\partial x}(t,x)\notag\\
    &+\frac{\partial }{\partial x}g(t, x, U^{\e}(t,x))+f(t, x, U^{\e}(t,x)),
                        \end{align}
for all $(t,x)\in[0,T]\times[0,1]$,  with Dirichlet boundary conditions
 ($U^{\e}(t,0)=U^{\e}(t,1)=0$) and initial condition $U^{\e}(0,x)=\eta(x)\in L^p([0,1]),p\ge2$;  $W$  denotes  the Brownian sheet  defined on a probability space $(\Omega,\FF,\{\FF_t\},\pp)$; the  coefficients  $f=f(t,x,r),g=g(t,x,r),\sigma=\sigma(t,x,r) $ are Borel functions of $(t,x,r)\in\rr_+\times[0,1]\times\rr$. See Section 2 for details.
 This family of semilinear equations contains both the stochastic Burgers' equation and the stochastic reaction-diffusion equation,
See Gy\"ongy \cite{Gy} for details.

 \vskip0.3cm

Intuitively, as the parameter $\e$ tends to zero, the solutions $U^{\e} $ of \eqref{SPDE} will tend to the solution of
\begin{align}\label{PDE}
       \frac{\partial U^{0}}{\partial t}(t,x)=&\frac{\partial^2 U^{0}}{\partial x^2}(t,x)+\frac{\partial }{\partial x}g(t, x, U^{0}(t,x))+f(t, x, U^{0}(t,x)),
\end{align}
for all $(t,x)\in[0,T]\times[0,1]$, with the Dirichlet's boundary conditions.

\vskip 0.3cm

It is  always interesting to  investigate deviations of $U^{\e} $ from the deterministic solution $U^0$, as $\e$ decreases to $0$, that is,  the asymptotic behavior of the trajectory,
$$
X^\e(t,x):=\frac{1}{\sqrt \e \lambda(\e)}\left(U^{\e} -U^{0}\right)(t,x),\quad (t,x)\in[0,T]\times[0,1],
$$
where $\lambda(\e)$ is some deviation scale which strongly influences the asymptotic behavior of $X^\e$.

\begin{itemize}
  \item[(1).]
 The case $\lambda(\e)=1/\sqrt\e$ provides some large deviations estimates.  Cardon-Weber \cite{CW} studied  the large deviations for the small noise limit of stochastic  semilinear PDEs by the exponential approximations. Very recently, Foondun and Setayeshgar \cite{FS} extended Cardon-Weber's result   to a less restrictive  case by using the weak convergence approach.

  \item[(2).]
  If $\lambda(\e)$ is identically equal to $1$, we are in the domain of the central limit theorem (CLT).
We will show that $(U^\e-U^0)/\sqrt\e$ converges  as $\e\downarrow 0$ to a  random field.
  \item[(3).]
 When the deviation scale satisfies
\beq \label{h}
 \lambda(\e)\to+\infty,\ \ \sqrt\e \lambda(\e)\to0\ \ \ \text{as}\ \e\to0,
 \nneq
 it is the moderate deviations, see \cite{DZ}. Throughout this paper, we assume \eqref{h} is in place.

\end{itemize}

 The moderate deviation principle (MDP) enables us to refine the estimates obtained through the  central limit theorem. It provides the asymptotic  behavior  for  $\pp(\|U^\e-U^0\| \ge \delta \sqrt \e h(\e))$ while the central limit theorem gives asymptotic bounds for
$\pp(\|U^\e-U^0\| \ge \delta \sqrt \e)$.

Like  large deviations, the moderate deviations arise in  the theory of statistical  inference quite naturally. The  moderate deviation principle   can provide us with the rate
  of convergence and a useful method for constructing asymptotic confidence intervals,  refer to the recent works   \cite{HS}, \cite{Gao} and the references therein. Usually, the quadratic form of the MDP's rate function
allows for the explicit minimization and in particular, it allows
to obtain an asymptotic evaluation for the exit time, see \cite{Klebaner1999}. Quite recently, the   study of the MDP estimates for  stochastic (partial) differential equation has been carried out as well, see \cite{BDG}, \cite{DXZZ}, \cite{GW}, \cite{LWYZ}, \cite{WZZ}, \cite{WZ} etc. Especially, Belfadli R. et al. \cite{BBM} proved a moderate deviation principle for the law of a stochastic Burgers equation and established some useful estimates toward a central limit theorem.

 \vskip0.3cm
 In this paper, we shall study the problems of  the central limit theorem and the  moderate deviation principle  for the stochastic semilinear SPDE \eqref{SPDE} which contains Burgers' equation and the stochastic reaction-diffusion equation. We generalize the moderate deviation result in \cite{BBM} and prove the central limit theorem.

 \vskip0.3cm

The rest of this paper is organized as follows. In Section 2, we give the framework of the the stochastic semilinear SPDEs, and   state the main results of this paper. In Section 3, we   prove some convergence results.  Section 4 is devoted to the proof of central limit theorem. In  Section 5, we prove the moderate deviation principle by using the weak convergence method.

\vskip0.3cm

Throughout the paper, $C(p)$ is a positive constant depending on the  parameter $p$, and $C$ is a constant depending on no specific parameter
(except $T$ and the Lipschitz constants), whose value may be different from line to line by convention.

\section{Framework and main results}

\subsection{Framework}

Let us give the framework taken from \cite{FS} and \cite{Gy}.

For  any $T>0$,    assume that the coefficients $f=f(t,x,r),g=g(t,x,r),\sigma=\sigma(t,x,r)$ in   \eqref{SPDE} are Borel functions of  $(t,x,r)\in[0,T]\times[0,1]\times\rr$  and there exist   positive constants $K,L$  satisfying  the following conditions:
\begin{itemize}
  \item[(H1)]  for all $(t,x,r)\in[0,T]\times[0,1]\times\rr$, it holds that
  $$|f(t,x,r)|\le K(1+|r|).
  $$
  \item[(H2)]the function $g$ is of the form $g(t,x,r)=g_1(t,x,r)+g_2(t,r)$, where $g_1$ and $g_2$ are Borel functions satisfying  that   $$
  |g_1(t,x,r)|\le K(1+|r|)\ \ \text{and} \ \ |g_2(t,r)|\le K(1+|r|^2).
  $$
  \item[(H3)] $\sigma$ is bounded and   for
  any $(t,x, p, q)\in [0,T]\times [0,1]\times \rr^2$,
  $$
  |\sigma(t,x,p)-\sigma(t,x,q)|\le L|p-q|.
  $$
  Furthermore, $f$ and $g$ are locally Lipschitz with linearly growing Lipschitz constant, i.e.,
  $$
  |f(t,x,p)-f(t,x,q)|\le L(1+|p|+|q|)|p-q|,
  $$
   $$
  |g(t,x,p)-g(t,x,q)|\le L(1+|p|+|q|)|p-q|.
  $$
\end{itemize}

\bdef[Mild Solution] A random field $U^{\e}=\{U^{\e}(t,x): t\in[0,T],x\in[0,1]\}$ is called a mild solution of \eqref{SPDE} with initial condition $\eta$ if $U^{\e}(t,x)$ is $\{\mathcal F_t\}$-measurable, $(t,x)\mapsto U^{\e}(t,x)$ is continuous a.s., and
\begin{align}\label{U e}
 U^{\e}(t,x) =&\int_0^1G_t(x,y)\eta(y)dy+\sqrt\e\int_0^t\int_0^1G_{t-s}(x,y)\sigma(s,U^{\e}(s))(y)W(dy,ds)\notag\\
 &-\int_0^t\int_0^1\partial_y G_{t-s}(x,y)g(s,U^{\e}(s))(y)dyds\notag\\
 &+\int_0^t\int_0^1 G_{t-s}(x,y)f(s,U^{\e}(s))(y)dyds.
 \end{align}
Here $G_t(\cdot,\cdot)$ is the Green kernel associated with the heat operator $ \partial/\partial t-\partial^2/\partial x^2$ with the  Dirichlet's  boundary conditions.
\ndef

  Gy\"ongy  \cite{Gy} proved the following result for the existence and uniqueness of the solution  to Eq.\eqref{SPDE}.
\bthm\label{thm solu}$($\cite[Theorem 2.1]{Gy}$)$. 
{\rm Under conditions (H1)-(H3), for any $\eta\in L^p([0,1]), p\ge2$, there exists a measurable functional
$$
\xi^\e: L^p([0,1])\times C([0,T]\times [0,1];\rr)\rightarrow C([0,T];L^p([0,1])),
$$
such that $U^{\e}=\xi^\e(\eta, \sqrt \e W)$ is the unique mild solution of \eqref{SPDE}.
}\nthm

Furthermore, from the proof of  \cite[Theorem 2.1]{Gy}, we know that   $\sup_{t\in[0,T]}\|U^\e(t,\cdot)\|_2$ is bounded in probability, i.e.,
\beq\label{eq bounded in prob}
\lim_{M\rightarrow \infty}\sup_{\e\in(0,1]}\pp\left(\sup_{t\in[0,T]}\left\|U^\e(t,\cdot)\right\|_2>M \right)=0.
\nneq
Particularly, taking $\e=0$ in \eqref{SPDE}, we know that the determinate equation \eqref{PDE} admits a unique solution  $U^0\in C([0,T];L^2([0,1]))$, given by
\begin{align}\label{U 0}
 U^{0}(t,x)=&\int_0^1G_t(x,y)\eta(y)dy-\int_0^t\int_0^1\partial_y G_{t-s}(x,y)g(U^{0})(s,y)dyds\notag\\
 &+\int_0^t\int_0^1 G_{t-s}(x,y)f(U^{0})(s,y)dyds,
 \end{align}
and
\begin{align}\label{U 02}
 \sup_{t\in[0,T]}\|U^{0}(t,\cdot)\|_{2}<\infty.
 \end{align}

\subsection{Main results}

To study the central limit theorem and moderate deviation principle, we furthermore suppose that
 \begin{itemize}
   \item[(H4)] the coefficients $f$ and $g$ are differentiable with respect to the last variable, and the derivatives $f'$ and $g'$ are also uniformly Lipschitz with respect to the last variable, more precisely, there exists a positive constant $K'$ such that
\beq\label{H4}
|f'(t,x,y)-f'(t,x,z)|\le K'|y-z|, \ \
|g'(t,x,y)-g'(t,x,z)|\le K'|y-z|
\nneq
for all $(t,x)\in[0,T]\times[0,1], y,z\in\rr$.
 \end{itemize}

Combined with the  growth condition (H3), we  conclude that
\beq\label{H4'}
|f'(t,x, r)|\le L(1+2|r|), \ \  \ \ |g'(t,x, r)|\le L(1+2|r|).
\nneq
\vskip0.3cm

Our first main  result is the following functional central limit theorem.
\bthm\label{CLT} 
{\rm Under conditions (H1)-(H4),   for any $T>0$, the process $(U^{\e}-U^0)/\e$  converges  to a random field $V$ in probability on $C([0,T];L^2([0,1]))$, determined by
\begin{align}\label{solu X}
V(t,x)=&\int_0^t\int_0^1G_{t-s}(x,y)f'(s,y, U^0(s,y))V(s,y)dyds\notag\\
&-\int_0^t\int_0^1\partial_y G_{t-s}(x,y)g'(s,y, U^0(s,y))V(s,y)dyds\notag\\
&+\int_0^t\int_0^1G_{t-s}(x,y)\sigma(s,y, U^0(s,y))W(dy,ds).
\end{align}

}
\nthm

\vskip0.3cm

In view of the assumption \eqref{h} and \eqref{H4}, by the large deviation principle for stochastic partial differential equation (see \cite{CW}), one can obtain that
 $V/\lambda(\e)$ obeys an LDP on $C([0,T];L^2([0,1]))$ with the speed $\lambda^2(\e)$ and with the good rate function:
\begin{equation}\label{rate}
I(f)=\inf\left\{\frac12\int_0^T\int_0^1|\dot{h}(t,x)|^2dtdx:\ h\in \HH,X^h=f\right\},
\end{equation}
where the function $X^h$ is the solution of the following  deterministic partial differential equation
\begin{align}\label{solu X h}
X^h(t,x)=&\int_0^t\int_0^1G_{t-s}(x,y)f'(s,y, U^0(s,y))X^h(s,y)dyds\notag\\
&-\int_0^t\int_0^1\partial_y G_{t-s}(x,y)g'(s,y, U^0(s,y))X^h(s,y)dyds\notag\\
&+\int_0^t\int_0^1G_{t-s}(x,y)\sigma(s,y, U^0(s,y))\dot h(s,y)dyds.
\end{align}

  Our second main result reads as follows:

 \bthm\label{thm MDP} {\rm Under conditions (H1)-(H4), $(U^{\e}-U^0)/(\sqrt\e\lambda(\e))$ obeys an LDP on $C([0,T];L^2([0,1]))$ with the speed $\lambda^2(\e)$ and with the rate function $I$ given by \eqref{rate}.

 }
 \nthm

\section{Some preliminary estimates}

\subsection{Some preliminary estimates}

The following estimates of Green function $G$ hold, (see Cardon-Weber \cite{CW}, Walsh \cite{Walsh}, Gy\"ongy \cite{Gy}). There exist positive constants $K,a,b,d$ such that for all $x,y\in[0,1]$ and $0\le s< t\le T$,
\begin{itemize}
  \item[(1)] \begin{equation}\label{formular1}|G_{t-s}(x,y)|\le K\frac{1}{ \sqrt{t-s}}\exp\left\{-a\frac{(x-y)^2}{t-s}\right\}, \end{equation}
  \item[(2)]\begin{equation}\label{formular2}\left|\frac{\partial }{\partial x}G_{t-s}(x,y)\right|\le K\frac{1}{ |t-s|^{3/2}}\exp\left\{-b\frac{(x-y)^2}{t-s}\right\}, \end{equation}
  \item[(3)]\begin{equation}\label{formular3}\left|\frac{\partial }{\partial t}G_{t-s}(x,y)\right|\le K\frac{1}{ |t-s|^{2}}\exp\left\{-d\frac{(x-y)^2}{t-s}\right\},\end{equation}
  \item[(4)]\begin{equation}\label{formular4}\sup_{t\in[0, T]}\int_0^t\int_0^1|G_u(x, z)-G_u(y, z)|^pdzdu\leq K|x-y|^{3-p}, \frac{3}{2}<p<3,\end{equation}
  \item[(5)]\begin{equation}\label{formular5}\sup_{x\in[0, 1]}\int_0^s\int_0^1|G_{t-u}(x, z)-G_{s-u}(x, z)|^pdzdu\leq K|t-s|^{(3-p)/2}, 1<p<3,\end{equation}
  \item[(6)]\begin{equation}\label{formular6}\sup_{x\in[0, 1]}\int_s^t\int_0^1|G_u(x, z)|^pdzdu\leq K|t-s|^{(3-p)/2}, 1<p<3.\end{equation}
\end{itemize}
For any $\tilde \alpha=\frac{\gamma-1}{2\gamma}$ with $\gamma\in (1,\infty), \alpha<\tilde \alpha$, there exists a constant $\tilde K(\alpha)$ such that for all $0<s<t<T, x,y\in[0,1]$,
\begin{itemize}
  \item[(7)] \begin{equation}\label{formular7}\int_0^T\int_0^1 |G_{t-u}(x,z)-G_{s-u}(y,z)|^2du dz\le \tilde K(\alpha)\rho((t,x),(s,y))^{2\alpha},\end{equation}
\end{itemize}
where $\rho$ is the Euclidean distance in $[0,T]\times[0,1]$.

For any  $v\in L^{\infty}([0,T], L^1([0,1])), t\in [0,T],x\in[0,1]$, let $J$ be a linear operator defined by
 $$
 J(v)(t,x):=\int_0^t\int_0^1 H(r,t;x,y) v(r,y)dydr,
 $$
with $H(t,s,x,y)=G_{t-s}(x,y), G_{t-s}^2(x,y)$ or  $\partial_yG_{t-s}(x,y)$.   Recall  the following regularity of $J(v)(t,x)$ from Gy\"ongy \cite{Gy}.

\blem$($\cite[Lemma 3.1]{Gy}$)$. \label{lem gyo}
{\rm
For any $\rho\in[1,\infty], q\in [1,\rho], K:=1+1/\rho-1/q$, it holds that $J$ is a bounded linear operator from $L^{\gamma}([0,T];L^q[0,1])$ into  $C([0,T];L^{\rho}[0,1])$ for any $\gamma>2K^{-1}$. Moreover, $J$ satisfies the following inequalities:
\begin{itemize}
  \item[(1).] For every $t\in[0, T]$ and $\gamma>2K^{-1}$,
  \begin{align*}
 \|J(v)(t,\cdot)\|_{\rho}\le&  C_1 \int_0^t(t-r)^{(1/2)K-1 }\|v(r,\cdot)\|_qdr\notag\\
 \le & C_2 t^{\frac K2-\frac1\gamma}\times\left(\int_0^t \|v(r,\cdot)\|_q^{\gamma}dr \right)^{\frac1\gamma}.
  \end{align*}
       \item[(2).] For every $0<\alpha< {K}/2$ and $\gamma>({K}/2-\alpha)^{-1}$, there exists a constant $C>0$ such that for all $0\le s\le t\le T$,
    $$
    \|J(v)(t,\cdot)-J(v)(s,\cdot)\|_{\rho}\le C (t-s)^{\alpha} \left(\int_0^t\|v(r,\cdot)\|_q^{\gamma} dr\right)^{\frac1\gamma}.
        $$
        \end{itemize}

}

\nlem

\subsection{The convergence of $U^\e$}

This section is concerned with the convergence of $U^\e$ to $U^0$ as $\e\to 0$.

For any  $M>0$, define the stopping time
$$
\tau^{M, \e}:=\inf\left\{t\ge 0;   \left\|U^{\e}(t,\cdot)\right\|_2 \ge M\right\}.
$$
 By \eqref{eq bounded in prob}, we know that
 \beq\label{eq stopping time}
\lim_{M\rightarrow \infty}\sup_{\e\in (0,1] } \pp\left(\tau^{M,\e}\le T\right)=0.
 \nneq

\bprop\label{Prop 2} 
{\rm Under conditions (H1)-(H3),  there exists some constant $C(M,L,\sigma, T)$ depending on $M,L,\sigma, T$ such that
\beq
\ee\left[\sup_{t\in[0,T\wedge\tau^{M,\e}]}\left\|U^\e(t,\cdot)-U^0(t,\cdot)\right\|_{2}^2\right]\le \e C(M,L,\sigma, T)\rightarrow 0,\ \ \text{as } \e\rightarrow0.
\nneq
}
\nprop
 \bprf
  Since
\begin{align}\label{Ue-U0}
 U^{\e}(t,x)-U^{0}(t,x)=& -\int_0^t\int_0^1\partial_y G_{t-s}(x,y)\left[g(U^{\e})-g(U^{0})\right](s,y)dyds\notag\\
 &+\int_0^t\int_0^1 G_{t-s}(x,y)\left[f(U^{\e})-f(U^{0})\right](s,y)dyds\notag\\
 &+ \sqrt{\e}\int_0^t\int_0^1G_{t-s}(x,y)\sigma(U^{\e})(s,y)W(dy,ds)\notag\\
 =:& I_1^\e+I_2^\e+I_3^\e.
 \end{align}
 By (H3) and Cauchy-Schwarz's inequality, we have that for any $(t,x)\in[0,T]\times[0,1]$,
\begin{align}\label{eq g}
 &\int_0^1\left|g(U^{\e})(s,y)-g(U^{0})(s,y)\right|dy\notag\\
 \le& L\int_0^1\left(1+\left|U^{\e}(s,y)\right|+\left|U^{0}(s,y)\right|\right)\cdot\left|U^{\e}(s,y)-U^{0}(s,y)\right|dy\notag\\
 \le & L\left(\int_0^1(1+|U^{\e}(s,y)|+|U^{0}(s,y)|)^2dy\cdot \int_0^1|U^{\e}(s,y)-U^{0}(s,y)|^2dy\right)^{\frac12}\notag\\
 \le & C(L)\left(\int_0^1(1+|U^{\e}(s,y)|^2+|U^{0}(s,y)|^2)dy\cdot \int_0^1|U^{\e}(s,y)-U^{0}(s,y)|^2dy\right)^{\frac12}.
 \end{align}
 Hence, applying    Lemma \ref{lem gyo} with $\rho=2, q=1$ and by  Cauchy-Schwarz's inequality,  we have
\begin{align*}
  &\|I_1^\e(s,\cdot)\|_2^2\notag\\
  \le &C\left(\int_0^s(s-r)^{-\frac34} \left\|g(U^{\e})-g(U^{0})(r,\cdot)\right\|_1dr \right)^{2}\notag\\
   \le & C(L)\left(\int_0^s(s-r)^{-\frac34}\left(1+\|U^{\e}(r,\cdot)\left\|_2^2+\|U^{0}(r,\cdot)\right\|_2^2 \right)^{\frac12} \cdot\|U^{\e}(r,\cdot)-U^{0}(r,\cdot)\|_2 dr\right)^{2}\notag\\
 \le & C(L) \int_0^s(s-r)^{-\frac34}\left(1+\|U^{\e}(r,\cdot)\|_2^2+\|U^{0}(r,\cdot)\|_2^2 \right)dr\cdot \int_0^s(s-r)^{-\frac34}\|U^{\e}(r,\cdot)-U^{0}(r,\cdot) \|_2^2dr.
 \end{align*}
 First taking the supremum of time over $[0,t\wedge \tau^{M, \e}]$, and then taking expectation, by the definition of $\tau^{M,\e}$, we
obtain that
 \begin{align}\label{I1}
 &\ee\left[\sup_{s\in[0,t\wedge \tau^{M, \e}]}\left\|I_1^\e(s,\cdot)\right\|_2^2\right]\notag\\
 \le & C(L) \ee\Bigg[\sup_{s\in[0,t\wedge \tau^{M, \e}]}\int_0^{s}(s-r)^{-\frac34}\left(1+\left\|U^{\e}(r,\cdot)\right\|_2^2+\left\|U^{0}(r,\cdot)\right\|_2^2 \right)dr\notag\\
 &\ \ \ \ \ \times \int_0^{s}(s-r)^{-\frac34}\left\|U^{\e}(r,\cdot)-U^{0}(r,\cdot) \right\|_2^2dr\Bigg]\notag\\
 \le & C(M,L) \ee\left[\sup_{s\in[0,t\wedge \tau^{M, \e}]}\int_0^{s}(s-r)^{-\frac34}\left\|U^{\e}(r,\cdot)-U^{0}(r,\cdot) \right\|_2^2dr\right]\notag\\
 \le& C(M,L)\int_0^{t}(t-r)^{-\frac34}\cdot \ee\left[\sup_{s\in[0,r\wedge \tau^{M, \e}]}\left\|U^{\e}(s,\cdot)-U^{0}(s\cdot) \right\|_2^2\right]dr,
  \end{align}
 where we have used the Fubini's theorem in the last inequality.

 For the term $I_2^\e$,  similarly to the proof of \eqref{I1}, we have
 \begin{align}\label{I2}
 &\ee\left[\sup_{s\in[0,t\wedge \tau^{M, \e}]}\left\|I_2^\e(s,\cdot)\right\|_2^2\right]\notag\\
   \le& C(M,L)\int_0^{t}(t-r)^{-\frac34}\cdot \ee\left[\sup_{s\in[0,r\wedge \tau^{M, \e}]}\left\|U^{\e}(s,\cdot)-U^{0}(s,\cdot) \right\|_2^2\right]dr.
  \end{align}

 For the term $I_3^\e$, first taking the supremum of time over $[0,t]$, and then taking expectation, by  Burkholder's inequality for stochastic integrals against Brownian sheets (see \cite{Ich}) and the boundedness of $\sigma$, we have
\begin{align}\label{I3}
 \ee\left[\sup_{s\in[0,t]}\left\|I_3^\e(s,\cdot)\right\|_2^2\right]
   \le&  \e \ee\left[\int_0^{t}\int_0^1G_{t-s}^2(x,y)\sigma^2(U^{\e})(s,y)dyds\right]\notag\\
 \le&\e C(M,\sigma).
  \end{align}

  Putting  \eqref{Ue-U0}, \eqref{I1}, \eqref{I2}, \eqref{I3} together,  we have
  \begin{align*}
  &\ee\left[\sup_{s\in[0,t\wedge \tau^{M, \e}]}\left\|U^{\e}(s,\cdot)-U^{0}(s,\cdot)\right \|_2^2\right]\\
   \le &\e C(M,\sigma)+C(M,L)\int_0^{t}(t-r)^{-\frac34} \ee\left[\sup_{s\in[0,r\wedge \tau^{M, \e}]}\left\|U^{\e}(s,\cdot)-U^{0}(s,\cdot) \right\|_2^2\right]dr.
  \end{align*}
  By  Gronwall's inequality (see \cite{YGD}), we know that
 there exists a constant $C(M,L,\sigma, T)$ such that
   \begin{align*}
   \ee\left[\sup_{s\in[0,t\wedge \tau^{M, \e}]}\|U^{\e}(s,\cdot)-U^{0}(s,\cdot) \|_2^2\right]\le \e C(M,L,\sigma, T).
   \end{align*}
The proof is complete.
 \nprf

\section{Proof of the Theorem \ref{CLT}}

 \bprf[Proof of Theorem \ref{CLT}]

 Let $$V^\e:=(U^\e-U^0)/\sqrt{\e}.$$
  We will prove that for any $\delta>0$,
  \beq\label{eq CLT}
  \lim_{\e\rightarrow0}\pp\left(\sup_{t\in[0,T]}\left\|V^\e(t,\cdot)-V(t,\cdot)\right\|_2\ge \delta \right)=0.
  \nneq

Recall  the stopping time defined by
$$
\tau^{M, \e}=\inf\left\{t\ge 0;   \left\|U^{\e}(t,\cdot)\right\|_2 \ge M\right\}.
$$
Since
\begin{align*}
&\pp\left(\sup_{t\in[0,T]}\left\|V^\e(t,\cdot)-V(t,\cdot)\right\|_2\ge \delta \right)\notag\\
\le &\pp\left(\sup_{t\in[0,T\wedge\tau^{M,\e}]}\left\|V^\e(t,\cdot)-V(t,\cdot)\right\|_2\ge \delta\right)+\pp(\tau^{M,\e}\le T),
\end{align*}
to prove \eqref{eq CLT}, by \eqref{eq stopping time}, it is sufficient to prove that for any  $M>0$ large enough
\begin{align}\label{eq CLT bound}
\lim_{\e\rightarrow0}\pp\left(\sup_{t\in[0,T\wedge\tau^{M,\e}]}\left\|V^\e(t,\cdot)-V(t,\cdot)\right\|_2\ge \delta\right)=0.
\end{align}

By the definition of $V^\e$ and $V$, we have
\begin{align}\label{eq Y e v 1}
&V^{\e}(t,x)-V(t,x)\notag\\
=&\int_0^t\int_0^1 G_{t-s}(x,y)\left[\frac{1}{\sqrt{\e}}\left(f(U^{\e})-f(U^0)\right)-f'(U^0)V\right](s,y)dyds\notag\\
&-\int_0^t\int_0^1\partial_yG_{t-s}(x,y)\left[\frac{1}{\sqrt{\e}}\left(g(U^{\e})-g(U^0)\right)-g'(U^0)V\right](s,y)dyds\notag\\
&+ \int_0^t\int_0^1G_{t-s}(x,y)\left[\sigma(U^{\e})  -\sigma( U^0)\right](s,y)W(dy,ds)\notag\\
=:&A_1^\e(t,x)+A_2^\e(t,x)+A_3^\e(t,x).
\end{align}

The first term $A_1^\e$ is further divided into two terms:
\begin{align}\label{eq A10}
&A_1^\e(t,x)\notag\\
=&\int_0^t\int_0^1 G_{t-s}(x,y)\left[\frac{1}{\sqrt{\e}}\left[f(U^{0}+\sqrt{\e} V^{\e})-f(U^0)\right]-f'(U^0)V^{\e}\right](s,y)dyds\notag\\
&+\int_0^t\int_0^1 G_{t-s}(x,y)\left[f'(U^0)(V^\e-V)\right](s,y)dyds\notag\\
=:&A_{11}^\e(t,x)+A_{12}^\e(t,x).
\end{align}
 By Taylor's formula, there exists a random field $\eta^\e$ taking values in $(0,1)$ such that
\begin{align*}
f(U^{0}+\sqrt{\e} V^{\e})-f(U^0)
= \sqrt{\e} f'\left( U^0+\sqrt{\e} \eta^{\e} V^{\e}\right)V^{\e}.
\end{align*}
Since $f'$ is  Lipschitz continuous, we have
\begin{align*}
\left|\frac{1}{\sqrt{\e}}\left[f(U^{0}+\sqrt{\e} V^{\e})-f(U^0)\right]-f'(U^0)V^{\e}\right|\le & \sqrt{\e} K' |V^\e|^2.
\end{align*}
Applying    Lemma \ref{lem gyo} for $A_{11}^\e$ with $\rho=2, q=1$, we have
$$
\|A_{11}^{\e}(s,\cdot)\|_2\le \sqrt{\e} K' \int_0^s(s-r)^{-\frac34}\|V^\e(r,\cdot)\|_2^2dr.
$$
 First taking the supremum of $s$ over $[0,t\wedge \tau^{M, \e}]$, and then taking expectation, by Proposition \ref{Prop 2}, we
 \begin{align}\label{eq A11}
\ee\left[\sup_{s\in[0,t\wedge \tau^{M, \e}]}\|A_{11}^\e(s,\cdot)\|_2\right]
\le& \sqrt{\e} K' \int_0^{t}(t-r)^{-\frac34}\ee\left[ \sup_{s\in [0,r\wedge \tau^{M, \e}]}\|V^\e(s,\cdot)\|_2^2\right]dr \notag\\
 \le & \sqrt{\e} C(K, K', L).
 \end{align}
 Applying    Lemma \ref{lem gyo} for  $A_{12}^\e$ with $\rho=2, q=1$, by Cauchy-Schwarz's inequality, we have
 \begin{align*}
 \|A_{12}^{\e}(s,\cdot)\|_2
 \le&  C\int_0^s (s-r)^{-\frac34}\left\|f'(U^0)(V^\e-V)(r,\cdot)\right\|_1dr\notag\\
 \le &  C\int_0^s (s-r)^{-\frac34}\left\|f'(U^0)(r,\cdot)\right\|_2\cdot\left\|(V^\e-V)(r,\cdot)\right\|_2dr.
  \end{align*}
First taking the supremum of $s$ over $[0,t\wedge \tau^{M, \e}]$, and then taking expectation, we have
 \begin{align}\label{eq A12}
 & \ee\left[\sup_{s\in[0,t\wedge \tau^{M, \e}]}\|A_{12}^\e(s,\cdot)\|_2\right]\notag\\
  \le & C\int_0^t (t-r)^{-\frac34}\left\|f'(U^0)(r,\cdot)\right\|_2\cdot\ee\left[\sup_{s\in [0,r\wedge \tau^{M, \e}]}\left\|(U^\e-U)(s,\cdot)\right\|_2\right]dr.
  \end{align}

 Putting \eqref{eq A11} and \eqref{eq A12}, we have
  \begin{align}\label{A1}
&\ee\left[\sup_{s\in[0,t\wedge \tau^{M, \e}]}\|A_{1}^{\e}(s,\cdot)\|_2\right]\notag\\
\le& \sqrt{\e} C_1+ C_2\int_0^t (t-r)^{-\frac34}\left\|f'(U^0)(r,\cdot)\right\|_2\cdot\ee\left[\sup_{s\in [0,r\wedge \tau^{M, \e}]}\left\|(V^\e-V)(s,\cdot)\right\|_2\right]dr.
 \end{align}

 For the term $A_2^\e$,  similarly to the proof of \eqref{A1}, we have
 \begin{align}\label{A2}
 &\ee\left[\sup_{s\in[0,t\wedge \tau^{M, \e}]}\|A_{2}^{\e}(s,\cdot)\|_2\right]\notag\\
 \le & \sqrt{\e} C_1+ C_2\int_0^t (t-r)^{-\frac34}\left\|g'(U^0)(r,\cdot)\right\|_2\cdot\ee\left[\sup_{s\in [0,r\wedge \tau^{M, \e}]}\left\|(V^\e-V)(s,\cdot)\right\|_2\right]dr.
   \end{align}

 For the term $A_3$, first taking the supremum of time over $[0,t\wedge \tau^{M, \e}]$, and then taking expectation, by  Burkholder's inequality for stochastic integrals against Brownian sheets (see \cite{Ich}), the Lipschitz continuity  of $\sigma$ and  Proposition \ref{Prop 2}, we have
 \begin{align}\label{A3}
 \ee\left[\sup_{s\in[0,t\wedge \tau^{M, \e}]}\|A_{3}^{\e}(s,\cdot)\|_2 \right]
 \le &   C\left(\int_0^t (t-r)^{\frac12}\ee\left[\sup_{s\in [0,r\wedge \tau^{M, \e}]}\left\|(U^\e-U^0)(s,\cdot)\right\|_2^2\right]dr\right)^{\frac12}\notag\\
 \le & \sqrt{\e} C(M,L,\sigma, T).
   \end{align}

Putting   \eqref{eq Y e v 1}, \eqref{A1}-\eqref{A3} together, we have
\begin{align*}
&\ee\left[\sup_{s\in [0,t\wedge \tau^{M, \e}]}\left\|(V^\e-V)(s,\cdot)\right\|_2\right]\\
\le &C_2\int_0^t (t-r)^{-\frac34}\left(\left\|f'(U^0)(r,\cdot)\right\|_2+\left\|g'(U^0)(r,\cdot)\right\|_2\right)\ee\left[\sup_{s\in [0,r\wedge \tau^{M, \e}]}\left\|(V^\e-V)(s,\cdot)\right\|_2\right]dr\\
&+\sqrt{\e} C_1(M,L,\sigma, T).
 \end{align*}
 By  Gronwall's inequalities (\cite[Theorem 1]{YGD}), \eqref{U 02} and \eqref{H4'}, we have
  \begin{align*}
\lim_{\e\rightarrow0}\ee\left[\sup_{t\in [0,T\wedge \tau^{M, \e}]}\left\|(V^\e-V)(t,\cdot)\right\|_2\right]=0,
    \end{align*}
    which implies \eqref{eq CLT bound}.

The proof is complete.

\nprf

 \section{Proof of the Theorem  \ref{thm MDP}}

\subsection{Weak convergence approach in LDP}\
First,  recall the  definition of  large deviation principle (c.f. \cite{DZ}).
Let $(\Omega,\mathcal{F},\mathbb{P})$ be a probability space with an increasing family $\{\FF_t\}_{0\le t\le T}$ of the sub-$\sigma$-fields of $\FF$ satisfying the usual conditions.
Let $\mathcal{E}$ be a Polish space with the Borel $\sigma$-field $\mathcal{B}(\mathcal{E})$.

    \bdef\label{def-Rate function}
         A function $I: \mathcal{E}\rightarrow[0,\infty]$ is called a rate function on
       $\mathcal{E}$,
       if for each $M<\infty$, the level set $\{x\in\mathcal{E}:I(x)\leq M\}$ is a compact subset of $\mathcal{E}$. A family of positive numbers $\{h(\e)\}_{\e>0}$ is called a speed function if $h(\e)\rightarrow +\infty$ as $\e\rightarrow 0$.
         \ndef

    \bdef\label{def LDP}
       $\{X^\e\}$   is  said to satisfy the large deviation principle on $\mathcal{E}$
       with rate function $I$ and with speed function $\{h(\e)\}_{\e>0}$, if the following two conditions
       hold:
       \begin{itemize}
         \item[$(a)$](Upper bound) For each closed subset $F$ of $\mathcal{E}$,
              $$
                \limsup_{\e\rightarrow 0}\frac1{h(\e)}\log\mathbb{P}(X^\e\in F)\leq- \inf_{x\in F}I(x);
              $$
         \item[$(b)$](Lower bound) For each open subset $G$ of $\mathcal{E}$,
              $$
                \liminf_{\e\rightarrow 0}\frac1{h(\e)}\log\mathbb{P}(X^\e\in G)\geq- \inf_{x\in G}I(x).
              $$
       \end{itemize}
    \ndef

% We shall apply the weak convergence approach to establish moderate deviation principle.

The  Cameron-Martin space  associated with the Brownian sheet $\{W(t,x);t\in[0,T],x\in[0,1]\}$ is given by
\begin{align}\label{eq CM}
\HH:=\left\{h(t,x)=\int_0^t\int_0^x \dot h(s,z)dzds; \dot h\in L^2([0,T]\times[0,1]) \right\}.
\end{align}
The space $\HH$ is a Hilbert space with inner product
$$\langle h_1,h_2\rangle_{\HH}:=\int_0^T\int_0^1 \dot h_1(s,z)\dot h_2(s,z)dzds.$$
The Hilbert space $\HH$ is endowed with   the norm $ \|h\|_{\HH}:=\left(\langle h,h\rangle_{\HH} \right)^{\frac12}$.

 Let $\AA$ denote  the class of real-valued $\{\FF_t\}$-predictable processes $\phi$ belonging to $\HH$ a.s., and let $$S_N:=\left\{h\in \HH; \|h\|_{\HH} \le N\right\}.$$
 The set $S_N$   endowed with the weak topology  is a Polish space.
Define $$\AA_N:=\left\{\phi\in \AA;\phi(\omega)\in S_N, \mathbb{P}\text{-a.s.}\right\}.$$

 Recall the following result from  \cite{BDM}.

\bthm\label{thm BD}$($\cite[Theorem 6]{BDM}$)$. {\rm For any $\e>0$, let $\Gamma^\e$ be a measurable mapping from $C([0,T];\rr)$ into $\EE$.
 Suppose that $\{\Gamma^\e\}_{\e>0}$ satisfies the following assumptions:
there  exists a measurable map $\Gamma^0:C([0,T];\rr)\longrightarrow \EE$ such that
\begin{itemize}
   \item[(a)] for any $N<+\infty$ and   family $\{h^\e;\e>0\}\subset \AA_N$ satisfying that $h^\e$ converge in distribution as $S_N$-valued random elements to $h$ as $\e\rightarrow 0$,
    $\Gamma^\e\left(W+\frac{1}{\sqrt\e}\int_0^{\cdot}\int_0^{\cdot}\dot h^\e(s,y)dyds\right)$ converges in distribution to  $\Gamma^0\left(\int_0^{\cdot}\int_0^{\cdot}\dot h(s,y)dyds\right)$ as $\e\rightarrow 0$;
   \item[(b)] for every $N<+\infty$, the set $
 \left\{\Gamma^0\left(\int_0^{\cdot}\int_0^{\cdot}\dot h(s,y)dyds\right); h\in S_N\right\}
  $   is a compact subset of $\EE$.
 \end{itemize}
Then the family $\{\Gamma^\e(W)\}_{\e>0}$ satisfies a large deviation principle in $\EE$ with the rate function $I$ given by
\beq\label{rate function}
I(g):=\inf_{\left\{h\in \HH;g=\Gamma^0\left(\int_0^{\cdot}\int_0^{\cdot}\dot h(s,y)dyds\right)\right\}}\frac12\|h\|_{\HH}^2, \  g\in\EE,
\nneq
with the convention $\inf \emptyset =\infty$.
 }\nthm

\subsection{The skeleton equation}
The purpose of this part is to study the skeleton equation, which will be used in the weak convergence approach.

Recall the  skeleton equation defined in \eqref{solu X h}. Using the same strategy in the proof of the existence and uniqueness for the solution to
Eq.\eqref{SPDE}, we know that
\begin{prop}\label{Prop skeleton}{\rm Under conditions (H1)-(H4), there exists a unique solution to Eq.\eqref{solu X h} satisfying that
\beq\label{eq Z h bound}
\sup_{h\in S_N}\sup_{t\in[0,T]}\left\|X^h(t,\cdot)\right\|_2<+\infty.
\nneq
}
\end{prop}

For any $h\in\HH$, set
\beq\label{Gamma 0} \Gamma^0\left(\int_0^\cdot\int_0^\cdot \dot h(s,y)dsdy\right):=X^h,
\nneq
where $X^h$ is the solution of   \eqref{solu X h}.

\bthm\label{thm Skeleton}{\rm  Under conditions (H1)-(H4), the mapping $h:S_N\rightarrow X^h\in C([0,T]; L^2([0,1]))$ is continuous with respect to the weak topology.
}
\nthm
\bprf Let  $\{h, (h_n)_{n\ge1}\}\subset S_N$ such that for any $g\in\HH$,
$$
\lim_{n\rightarrow\infty}\langle h_n-h, g \rangle_{\HH}=0.
$$
 We need to  prove that
\beq\label{eq convergence}
\lim_{n\rightarrow \infty}\sup_{t\in[0,T]}\left\|X^{h_n}(t, \cdot)-X^h(t,\cdot)\right\|_2=0.
\nneq

Notice that
\begin{align}\label{eq ske  1}
&X^{h_n}(t,x)-X^h(t,x)\notag\\
=&\int_0^t\int_0^1G_{t-s}(x,y)f'(s,y, U^0(s,y))\left(X^{h_n}(s,y)-X^{h}(s,y)\right)dyds\notag\\
&-\int_0^t\int_0^1\partial_yG_{t-s}(x,y)g'(s,y, U^0(s,y))\left(X^{h_n}(s,y)-X^{h}(s,y)\right)dyds\notag\\
&+\int_0^t\int_0^1G_{t-s}(x,y)\sigma(s,y, U^0(s,y))\left(\dot h_n(s,y)-\dot h(s,y)\right)dyds\notag\\
=:&I_1^n(t,x)+I_2^n(t,x)+I_3^n(t,x).
\end{align}
 By \eqref{U 02}, \eqref{H4'} and \eqref{eq Z h bound},  we have
\begin{align*}
&\sup_{s\in[0,T]}\int_0^1\left|f'(s,y,U^0(s,y))\left(X^{h_n}(s,y)-X^{h}(s,y)\right)\right|dy\\
\le&\sup_{s\in[0,T]}\left(\int_0^1\left|f'(s,y,U^0(s,y)) \right|^2dy\right)^{\frac12}\cdot
\left(\int_0^1 \left|X^{h_n}(s,y)-X^{h}(s,y) \right|^2dy\right)^{\frac12}\\
\le& \sup_{s\in[0,T]}\left(\int_0^1\left( L(1+2|U^0(s,y)|)\right)^2 dy\right)^{\frac12}\cdot
\left(\int_0^1 \left(2|X^{h_n}(s,y)|^2+2|X^{h}(s,y)|^2\right)dy\right)^{\frac12}\\
< & \infty.
\end{align*}
Thus, the function  $(s,y)\mapsto f'(s,y,U^0(s,y))\left(X^{h_n}(s,y)-X^{h}(s,y)\right)$   belongs  to $L^{\infty}([0,T]; L^1([0,1]))$. Applying Lemma \ref{lem gyo}  with $\rho=2, q=1$ and by the H\"older inequality, we have
\begin{align}\label{eq ske 11}
\|I_1^n(t,\cdot)\|_2\le & C\int_0^t(t-s)^{-\frac34}\cdot\|X^{h_n}(s,\cdot)-X^{h}(s,\cdot)\|_1ds\notag\\
 \le & C\int_0^t(t-s)^{-\frac34}\cdot\|X^{h_n}(s,\cdot)-X^{h}(s,\cdot)\|_2ds.
\end{align}
Similarly, we obtain that
\beq\label{eq ske 12}
\|I_2^n(t,\cdot)\|_2\le C\int_0^t(t-s)^{-\frac34}\cdot\|X^{h_n}(s,\cdot)-X^{h}(s,\cdot)\|_2ds.
\nneq
 Since $\sigma$ is bounded, for any fixed  $(t,x)\in [0,T]\times [0,1]$, the function $G_{t-s}(x,y)\sigma(s,y, U^0(s,y)):(s,y)\in[0,t]\times[0,1]\rightarrow \rr$ belongs to  $L^2([0,T]\times[0,1];\rr)$.  As $\dot h_n\rightarrow \dot h$ weakly in $L^2([0,T]\times[0,1];\rr)$, it holds that
 \begin{align}\label{eq ske 131}
 I^n_3(t,x)=\int_0^t\int_0^1G_{t-s}(x,y)\sigma(s,y, U^0(s,y))\left(\dot h_n(s,y)-\dot h(s,y)\right)dyds\rightarrow 0.
\end{align}

 For any $0\le s\leq t \le T$, applying formulars \eqref{formular1}, \eqref{formular5} and \eqref{formular6}, by the boundness of $\sigma$ and H\"{o}lder's inequality, we obtain that
\begin{align*}
&|I_3^n(t,x)-I_3^n(s,x)|\notag\\
\le& \left|\int_0^s \int_0^1(G_{t-u}(x, y)-G_{s-u}(x, y)) \sigma(u,y, U^0(u,y))(\dot h_n(u,y)-\dot h(u,y))dudy \right|\\
&+\left|\int_s^t \int_0^1G_{t-u}(x, y) \sigma(u,y, U^0(u,y))(\dot h_n(u,y)-\dot h(u,y))dudy \right|\notag\\
\le& C(N, \sigma)\left(\left(\int_0^s \int_0^1|G_{t-u}(x, y)-G_{s-u}(x, y)|^2dudy\right)^{\frac{1}{2}}+\left(\int_s^t \int_0^1G_{t-u}(x, y)^2dudy \right)^{\frac{1}{2}}\right)\notag\\
\le& C(N, \sigma)(t-s)^{\frac{1}{4}}.
\end{align*}
Then
\begin{align*}
\|I_3^n(t,\cdot)-I_3^n(s,\cdot)\|_2 \le C(N, \sigma)(t-s)^{\frac{1}{4}}.
\end{align*}
Particularly, taking $s=0$, we obtain that
\beq\label{eq I21}
\|I_3^n(t,\cdot)\|_2\le C(N, \sigma)t^{\frac{1}{4}}.
\nneq
 Hence, the functions $\{I_3^n\}_{n\ge1}$ are uniformly bounded and equi-continuous in $C([0,T]; L^2([0,1]))$.    According to Arz\'ela-Ascoli theorem,  the functions  $t\mapsto\{ I_3^n(t,\cdot)\}_{n\ge1}$ is pre-compact in $C([0,T]; L^2([0,1]))$. Thus, by \eqref{eq ske 131}, we obtain that
\beq\label{eq I2}
\lim_{n\rightarrow\infty}\sup_{t\in[0,T]}\|I^n_3(t,\cdot)\|_2=0.
 \nneq

  Set $\zeta^n(t)=\sup_{0\le s\le t}\left\|X^{h_n}(s, \cdot)-X^h(s,\cdot)\right\|_2$. By \eqref{eq ske  1}-\eqref{eq ske 12}, we have
   \begin{align*}
   \zeta^n(t)\le C\int_0^t(t-s)^{-\frac34} \zeta^n(s) ds+\sup_{s\in[0,t]}\|I^n_3(s,\cdot)\|_2.
\end{align*}
  Hence, by a generalized Gronwall  lemma (eg. \cite[Theorem 1]{YGD}), we have
  $$
  \zeta^n(t)\le C\sup_{s\in[0,t]}\|I^n_3(s,\cdot)\|_2,
  $$
   which, together with \eqref{eq I2}, implies the desired estimate \eqref{eq convergence}.

The proof is complete.
\nprf

\subsection{The proof of Theorem \ref{thm MDP}}

Recall that  $X^\e=(U^\e-U^0)/(\sqrt\e \lambda(\e))$. By  \eqref{U e} and \eqref{U 0}, we know that
\begin{align}\label{eq X e}
&X^{\e}(t,x)\notag\\
=&\frac{1}{\lambda(\e)}\int_0^t\int_0^1G_{t-s}(x,y)\sigma(s,y, U^0(s,y)+\sqrt{\e}\lambda(\e)X^{\e}(s,y))W(dy,ds)\notag\\
&+\frac{1}{\sqrt\e\lambda(\e)} \int_0^t\int_0^1G_{t-s}(x,y)\left(f(s,y, U^0(s,y)+\sqrt{\e}\lambda(\e)X^{\e}(s,y))-f(s,y, U^0(s,y))\right)dyds\notag\\
&- \frac{1}{\sqrt\e\lambda(\e)} \int_0^t\int_0^1\partial_yG_{t-s}(x,y)\left(g(s,y, U^0(s,y)+\sqrt{\e}\lambda(\e)X^{\e}(s,y))-g(s,y, U^0(s,y))\right)dyds.
\end{align}
 This equation admits a unique strong solution  \beq\label{gamma e}
   X^\e:=\Gamma^{\e}(W),
   \nneq
  where $\Gamma^\e$ stands for the solution
functional from $C([0,T]\times [0,1];\rr)$ into $C([0,T]\times[0,1];\rr)$.

The following lemma  is a direct consequence of Girsanov's theorem, refer to \cite[Theorem 3.2]{FS}.
\blem\label{lem control} {\rm For every fixed  $N\in\mathbb{N}$,    let $v\in \mathcal{A}_N$ and $\Gamma^\e$ be given by \eqref{gamma e}. Then  $X^{\e,v}:=\Gamma^\e\left(W+\lambda(\e)v\right) \in C([0,T];L^2([0,1]))$    solves the following equation:
\begin{align}\label{eq X e v}
&X^{\e,v}(t,x)\notag\\
=&\frac{1}{\lambda(\e)}\int_0^t\int_0^1G_{t-s}(x,y)\sigma(s,y, U^0(s,y)+\sqrt{\e}\lambda(\e)X^{\e,v}(s,y))W(dy,ds)\notag\\
&+\frac{1}{\sqrt\e\lambda(\e)} \int_0^t\int_0^1G_{t-s}(x,y)\left(f(s,y, U^0(s,y)+\sqrt{\e}\lambda(\e)X^{\e,v}(s,y))-f(s,y, U^0(s,y))\right)dyds\notag\\
&-\frac{1}{\sqrt\e\lambda(\e)} \int_0^t\int_0^1\partial_yG_{t-s}(x,y)\left(g(s,y, U^0(s,y)+\sqrt{\e}\lambda(\e)X^{\e,v}(s,y))-g(s,y, U^0(s,y))\right)dyds\notag\\
&+\int_0^t\int_0^1G_{t-s}(x,y)\sigma(s,y, U^0(s,y)+\sqrt{\e}\lambda(\e)X^{\e,v}(s,y))\dot v(s,y)dyds.
\end{align}
Furthermore,
\beq\label{eq Z e estimate}
\lim_{M\rightarrow\infty}\sup_{\e\le 1}\sup_{v\in\AA_N}\pp\left(\sup_{t\in[0,T]}\|X^{\e,v}(t,\cdot)\|_2^2\geq M\right)=0.
\nneq
}\nlem
\vskip0.3cm

We are now ready to state our main result. Recall the mapping $\Gamma_0$ given by \eqref{Gamma 0}. For any $g\in C([0,T];L^2([0,1]))$, let
\beq\label{eq rate}
I(g):=\inf_{\left\{h\in \HH; g=\Gamma^0\left(\int_0^{\cdot}\int_0^{\cdot}\dot h(s,y)dyds \right)\right\}}\frac{1}{2}\|h\|_{\HH}^2.
\nneq

\bprf  According  to Theorem \ref{thm BD},  we need to prove that the following two conditions  are fulfilled:
 \begin{itemize}
   \item[(a)] the set $\{X^h;h\in S_N\}$ is a compact set of $C([0,T];L^2([0,1]))$, where $X^h$ is the solution of Eq.\eqref{solu X h}.
   \item[(b)] for any family $\{v^\e;\e>0\}\subset\AA_N$ which converges in distribution as $\e\rightarrow 0$ to $v\in \AA_{N}$, as $S_N$-valued random variables, we have
       $$
       \lim_{\e\rightarrow 0}X^{\e,v^{\e}}=X^v\ \text{ in distribution,}
       $$
 as $C([0,T];L^2([0,1]))$-valued random variables, where $X^v$ denotes the solution of Eq.\eqref{solu X h} corresponding to the $S_N$-valued random variable $v$ (instead of a deterministic function $h$).
 \end{itemize}
  Condition (a)  follows from the continuity of the mapping $h:S_N\rightarrow X^h\in C([0,T];L^2([0,1]))$, which has been   established in Theorem \ref{thm Skeleton}. The verification of condition (b) will be given by Proposition \ref{prop convergence} below.

  The proof is complete.
\nprf
\vskip0.3cm

\bprop\label{prop convergence}{\rm Assume (H1)-(H4). For every fixed  $N\in\mathbb{N}$, let $v^\e,v\in \mathcal{A}_N$ be such that $v^\e$ converges in
distribution to $v$ as $\e\longrightarrow0$. Then
$$
\Gamma^\e\left(W+\lambda(\e)v^\e \right)\ \ \text{converges in
distribution to } \ \Gamma^0\left(v\right),
$$
in $C([0,T];L^2([0,1]))$ as $\e\longrightarrow0$.

}
\nprop
\bprf
By   Skorokhod representation theorem, there exist a probability basis $(\bar\Omega,\bar\FF,(\bar\FF_t),\bar\pp)$, and, on this basis, a sequence of independent Brownian sheets $\bar W=(\bar W_k)_{k\ge1}$ and also a family of $\bar \FF_t$-predictable processes $\{\bar v^\e;\e>0\}, \bar v$ belonging to $L^2(\bar\Omega\times [0,T];\HH)$ taking values on $S_N$, $\bar\pp$-a.s., such that the joint law of $(v^\e,v, W)$ under $\pp$ coincides with that of  $(\bar v^\e,\bar v, \bar W)$ under $\bar\pp$ and
$$
\lim_{\e\rightarrow0}\langle\bar v^\e-\bar v, g \rangle_{\HH}=0, \ \ \forall g\in\HH, \ \bar \pp\text{-a.s.}.
$$

Let $\bar X^{\e,\bar v^\e}$ be the solution to a similar equation as \eqref{eq X e v} replacing $v$  by $\bar v^\e$ and $W$ by $\bar W$. Thus, to prove this proposition, it is sufficient to prove that
\beq\label{eq b conv}
\lim_{\e\rightarrow0}  \sup_{0\le t\le T}\left\|\bar X^{\e,\bar v^\e}-\bar X^{\bar v} \right\|_2=0, \ \ \ \ \text{in probability}.
\nneq

From now on, we drop the bars in the notation for the sake of simplicity, and we denote
$$
Y^{\e,v^\e,v}:=X^{\e,v^\e}-X^v.
$$
Notice that
\begin{align}\label{eq Y e v}
&X^{\e,v^\e}(t,x)-X^{v}(t,x)\notag\\
=&\frac{1}{\lambda(\e)}\int_0^t\int_0^1G_{t-s}(x,y)\sigma(s,y, U^0(s,y)+\sqrt{\e}\lambda(\e)X^{\e,v^\e}(s,y))W(dy,ds)\notag\\
&+\int_0^t\int_0^1 G_{t-s}(x,y)\Bigg[\frac{1}{\sqrt\e\lambda(\e)}\left(f(s,y, U^0(s,y)+\sqrt{\e}\lambda(\e)X^{\e,v^\e}(s,y))-f(s,y, U^0(s,y))\right)\notag\\
&\ \ \ \ \ \ \ \ \ \ \ \ \ \ \ \ \ \ \ \ \ \ \ \ \ \ \ \ - f'(s,y, U^0(s,y))X^{v}(t,x)\Bigg]dyds\notag\\
&- \int_0^t\int_0^1\partial_yG_{t-s}(x,y)\Bigg[\frac{1}{\sqrt\e\lambda(\e)}\left(g(s,y, U^0(s,y)+\sqrt{\e}\lambda(\e)X^{\e,v^\e}(s,y))-g(s,y, U^0(s,y))\right)\notag\\
&\ \ \ \ \ \ \ \ \ \ \ \ \ \ \ \ \ \ \ \ \ \ \ \ \ \ \ \ - g'(s,y, U^0(s,y))X^{v}(t,x)\Bigg]dyds\notag\\
&+\int_0^t\int_0^1G_{t-s}(x,y)\Bigg[\sigma(s,y, U^0(s,y)+\sqrt{\e}\lambda(\e)X^{\e,v^\e}(s,y)) \dot v^\e(s,y)\notag\\
&\ \ \ \ \ \ \ \ \ \ \ \ \ \ \ \ \ \ \ \ \ \ \ \ \ \ \ \ - \sigma(s,y, U^0(s,y))\dot v(s,y) \Bigg]dyds\notag\\
=:&A_1^\e(t,x)+A_2^\e(t,x)+A_3^\e(t,x)+A_4^\e(t,x).
\end{align}
\vskip0.3cm
{\bf Term $A_1^\e$.} Since $\sigma$ is bounded,  by Burkholder's inequality for stochastic integrals against Brownian sheets (see  \cite{Ich}), we have
\begin{align}\label{eq Y e v1}
 \ee\left[\sup_{t\in[0,T]}\|A_2^\e(t,\cdot)\|_2^2\right]\leq  \frac{1}{\lambda^2(\e)}C(\sigma,T)\int_0^T\int_0^1G_{t-s}^2(x,y) dyds \longrightarrow0, \ \ \text{as } \e\rightarrow0.
\end{align}
\vskip0.3cm
{\bf Terms $A_2^\e$ and  $A_3^\e$.}  Notice that
\begin{align}\label{eq Y e v20}
&A_2^\e(t,x)\notag\\
=&\int_0^t\int_0^1 G_{t-s}(x,y)\Bigg[\frac{1}{\sqrt\e\lambda(\e)}\left(f(s,y, U^0(s,y)+\sqrt{\e}\lambda(\e)X^{\e,v^\e}(s,y))-f(s,y, U^0(s,y))\right)\notag\\
&\ \ \ \ \ \ \ \ \ \ \ \ \ \ \ \ \ \ \ \ \ \ \ \ \ \ \ \ - f'(s,y, U^0(s,y))X^{\e,v^\e}(s,y)\Bigg]dyds\notag\\
&+\int_0^t\int_0^1 G_{t-s}(x,y)\Bigg[ f'(s,y, U^0(s,y))(X^{\e,v^{\e}}(s,y)-X^{v}(s,y))\Bigg]dyds\notag\\
=:&A_{21}^\e(t,x)+A_{22}^\e(t,x).
\end{align}

 By Taylor's formula, there exists a random field $\eta^\e(s,y)$ taking values in $(0,1)$ such that
\begin{align*}
&f\left(s,y, U^0(s,y)+\sqrt{\e}\lambda(\e)X^{\e,v^\e}(s,y)\right)-f\left(s,y, U^0(s,y)\right)\\
=&\sqrt\e\lambda(\e) f'\left(s,y, U^0(s,y)+ \sqrt{\e}\lambda(\e) \eta^{\e}(s,y) X^{\e,v^\e}(s,y)\right)X^{\e,v^\e}(s,y).
\end{align*}
Since $f'$ is also Lipschitz continuous, we have
\begin{align*}
&\left|f'(s,y, U^0(s,y)+ \sqrt{\e}\lambda(\e) \eta^{\e}(s,y) X^{\e,v^\e}(s,y))X^{\e,v^\e}(s,y)-f'(s,y, U^0(s,y))X^{\e,v^\e}(s,y)\right|\\
\le&\ K \sqrt\e\lambda(\e) \left|X^{\e,v^\e}(s,y) \right|^2.
\end{align*}

Define the stopping time
$$
\tau^{M, \e}:=\inf\left\{t\ge 0;   \|X^{\e,v^\e}(t,\cdot)\|_2\vee \|X^{v}(t,\cdot)\|_2\ge M\right\},
$$
where $M$ is some constant large enough.

 Applying Lemma \ref{lem gyo}  for $A_{21}^\e$ with $\rho=2, q=1$,   we have that for all $t\in [0,T]$,
\begin{align*}
 \|A_{21}^\e(t,\cdot)\|_2
\le  \sqrt\e\lambda(\e)C(K) \int_0^t(t-s)^{-\frac34}\cdot\left\|(X^{\e,v^\e}(s,\cdot))^2\right\|_1ds.
\end{align*}
First taking the supremum of time over $[0,t\wedge \tau^{M, \e}]$, and then taking expectation, we  obtain that
\begin{align}\label{eq Y e v201}
&\ee \left[\sup_{s\in[0,t\wedge \tau^{M, \e}]}\|A_{21}^\e(s,\cdot)\|_2 \right]\notag\\
\le & \sqrt\e\lambda(\e)C(K) \int_0^t(t-s)^{-\frac34}\ee\left[\sup_{0\le r\le s\wedge \tau^{M, \e}}\left\|X^{\e,v^\e}(r,\cdot)\right\|_2^2\right]ds\notag\\
\le & \sqrt\e\lambda(\e)C(K, T,M).
\end{align}

For the  term $A_{22}^\e$, similarly to  the proof of \eqref{eq ske 11},   we have
\begin{align}\label{eq Y e v202}
&\ee \left[\sup_{s\in[0,t\wedge \tau^{M, \e}]}\|A_{22}^\e(s,\cdot)\|_2 \right]\notag\\
\le& C\int_0^t(t-s)^{-\frac34}\cdot\ee\left[\sup_{0\le r\le s\wedge \tau^{M, \e}}\left\|X^{\e,v^\e}(r,\cdot)-X^{v}(r,\cdot)\right\|_2\right]ds.
\end{align}

Putting \eqref{eq Y e v201} and \eqref{eq Y e v202} together, we have
\begin{align}\label{eq Y e v2}
& \ee\left[\sup_{s\in[0,t\wedge \tau^{M, \e}]} \|A_{2}^\e(s,\cdot)\|_2\right]\notag\\
\le &  \sqrt\e\lambda(\e)C(K, T,M)+ C\int_0^t(t-s)^{-\frac34}\cdot\ee\left[\sup_{0\le r\le s\wedge \tau^{M, \e}}\left\|X^{\e,v^\e}(r,\cdot)-X^{v}(r,\cdot)\right\|_2\right]ds.
\end{align}

Similarly to  the proof of \eqref{eq Y e v2}, we obtain the following estimate for $A_3^\e$:
\begin{align}\label{eq Y e v3}
& \ee\left[\sup_{s\in[0,t\wedge \tau^{M, \e}]} \|A_{3}^\e(s,\cdot)\|_2\right]\notag\\
\le &  \sqrt\e\lambda(\e)C(K, T,M)+ C\int_0^t(t-s)^{-\frac34}\cdot\ee\left[\sup_{0\le r\le s\wedge \tau^{M, \e}}\left\|X^{\e,v^\e}(r,\cdot)-X^{v}(r,\cdot)\right\|_2\right]ds.
\end{align}

{\bf Term $A_4^\e$.}  Notice that
\begin{align}\label{eq Y e v40}
&A_4^\e(t,x)\notag\\
=&\int_0^t\int_0^1G_{t-s}(x,y)\left[\sigma(s,y, U^0(s,y)+\sqrt{\e}\lambda(\e)X^{\e,v^\e}(s,y)) ( v^\e(s,y)-v(s,y)\right]dyds\notag\\
& +\int_0^t\int_0^1G_{t-s}(x,y)\left\{\left[\sigma(s,y, U^0(s,y)+\sqrt{\e}\lambda(\e)X^{\e,v^\e}(s,y))-   \sigma(s,y, U^0(s,y))\right]\dot v(s,y) \right\}dyds\notag\\
=:&A_{41}^\e(t,x)+A_{42}^\e(t,x).
\end{align}

Using the  argument  as that in the proof of  \eqref{eq I2},  we obtain that
 \begin{align}\label{eq Y e v41}
 \lim_{\e\rightarrow0}\ee\left[\sup_{t\in[0,T]}\left\|A_{41}^\e(t,\cdot)\right\|_2\right]= 0.
\end{align}
By the Lipschitz continuity of $\sigma$, we know that
\begin{align*}
  | A_{42}^\e(t,x)|\le \sqrt{\e}\lambda(\e)K \int_0^t\int_0^1G_{t-s}(x,y)|X^{\e,v^\e}(s,y)|\cdot| \dot v(s,y)| dyds.
\end{align*}
 Applying Lemma \ref{lem gyo} for $A_{42}^\e$  with $\rho=2, q=1$ and by H\"older's inequality,  we obtain that for all $t\in [0,T]$,
\begin{align*}
\|A_{42}^\e(t,\cdot)\|_2
\le & \sqrt\e\lambda(\e)C(K) \int_0^t(t-r)^{-\frac34}  \|X^{\e,v^\e}(r,\cdot)    \dot v(r,\cdot)\|_1 dr\notag\\
\le & \sqrt\e\lambda(\e)C(K) \int_0^t(t-r)^{-\frac34}  \|X^{\e,v^\e}(r,\cdot)\|_2  \cdot  \|\dot v(r,\cdot)\|_2 dr.
\end{align*}

First taking the supremum of $t$ over $[0,T\wedge \tau^{M, \e}]$, and then taking expectation, we  obtain that
\begin{align}\label{eq Y e v42}
&\ee\left[\sup_{t\in[0,T]} \|A_{42}^\e(t\wedge \tau^{M, \e},\cdot)\|_2 \right]\notag\\
\le & \sqrt\e\lambda(\e)C(K) \int_0^t(t-r)^{-\frac34}\ee\left[\sup_{0\le r\le s\wedge \tau^{M, \e}}\left\|X^{\e,v^\e}(r,\cdot)\right\|_2 \cdot \|\dot v(r,\cdot)\|_2\right]dr\notag\\
\le & \sqrt\e\lambda(\e)C(K, T, N,M).
\end{align}

According to \eqref{eq Y e v40}-\eqref{eq Y e v42}, we have that
\begin{align}\label{eq Y e v4}
 \lim_{\e\rightarrow0} \ee\left[\sup_{t\in[0,T]}\left\|A_{4}^\e(t\wedge \tau^{M, \e},\cdot)\right\|_2\right]= 0.
\end{align}

Putting   \eqref{eq Y e v}, \eqref{eq Y e v1}, \eqref{eq Y e v2}, \eqref{eq Y e v3} and \eqref{eq Y e v4} together   and by  Gronwall's inequality (\cite[Theorem 1]{YGD}), we have
  \begin{align}\label{eq Y e v5}
\lim_{\e\rightarrow0}\ee\left[\sup_{t\in [0,T]}\left\|X^{\e,v^\e}(s\wedge \tau^{M, \e},\cdot)-X^{v}(t\wedge \tau^{M, \e},\cdot)\right\|_2\right]=0.
    \end{align}
 By Chebychev's inequality, we have
 $$
 \lim_{\e\rightarrow 0} \sup_{t\in [0,T]}\left\|X^{\e,v^\e}(s\wedge \tau^{M, \e},\cdot)-X^{v}(t\wedge \tau^{M, \e},\cdot)\right\|_2=0, \ \ \text{in probability}.
$$
Letting $M\rightarrow\infty$, by (\ref{eq Z e estimate}) we get the desired result \eqref{eq b conv}.

The proof is complete.
\nprf
%The following generalized Gronwall inequality is studied by Ye {\it et al.}  in \cite{YGD}.
%\bthm[\cite{YGD}]{\rm Suppose $\beta>0$, $a(t)$ is a nonnegative function locally integrable on $0\le t<T$ (some $T\le +\infty$) and $g(t)$ is a nonnegative, nonnegative, nondecreasing continuous function defined on $0\le t<T$, $g(t)\le M$ (constant), and suppose $u(t)$ is nonnegative and locally integrable on $0\le t<T$ with
%$$
%u(t)\le a(t)+g(t)\int_0^t (t-s)^{\beta-1}u(s)ds
%$$
%on this interval. Then
%$$
%u(t)\le a(t)+\int_0^t\left[\sum_{n=1}^{\infty}\frac{(g(t)\Gamma(\beta))^n}{\Gamma(n\beta)}(t-s)^{n\beta-1}a(s) \right]ds, \ \ \ 0\le t<T.
%$$

%}\nthm

\section*{Acknowledgement}
We would like to express our appreciation for Belfadli R., Boulanba L. and Mellouk M.. They told us their work about the stochastic Burgers equation and pointed out a  mistake in our paper.

\vskip0.3cm
%\noindent{\bf Acknowledgements:} R. Wang was supported by Natural Science Foundation of China %(11301498, 11431014, 11671076).


\begin{thebibliography}{abc}
%\bibitem{BF} Baier, U., Freidlin, M.I.: Theorems on large deviations and stability under random perturbations. Dokl. Akad. Nauk USSR 235, 253-256 (1977)

% \bibitem{BMS}Bally, V., Millet, A., Sanz-Sol\'{e}, M.:  Approximation and support theorem in H\"{o}lder normfor parabolic stochastic partial differential equations. Ann. Probab. 23, 178-222 (1995)

\bibitem{BBM} Belfadli R., Boulanba L., Mellouk M.:  Moderate deviations for a stochastic Burgers equation.   	arXiv:1807.09117.

\bibitem{BDG} Budhiraja A., Dupuis P., Ganguly A.:  Moderate deviation principle for stochastic differential equations with jumps.   Ann. Probab. 44, 1723-1775 (2016).

\bibitem{BDM}Budhiraja  A., Dupuis  P., Maroulas  V.: Large deviations for infinite dimensional stochastic dynamical systems.   Ann.  Probab. 36, 1390-1420 (2008).



\bibitem{CW} Cardon-Weber C.: Large deviations for Burgers' type SPDE.  Stochastic Process. Appl. 84, 53-70 (1999).


%\bibitem{DPZ} Da Prato, G., Zabczyk, J.: Stochastic equations in infinite dimensions. Cambridge Univ. Press, Cambridge (1992)

%\bibitem{DeA} De Acosta, A.: Moderate deviations and associated Laplace approximations for sums of independent random vectors.
%Trans. Amer. Math. Soc. 329, 357-375 (1992)

\bibitem{DZ} Dembo  A., Zeitouni  O.:  Large deviations techniques and applications. Second
edition. Applications of Mathematics 38, Springer-Verlag (1998).

 \bibitem{DXZZ} Dong  Z., Xiong  J., Zhai  J., Zhang  T.: A moderate deviation principle for 2-D stochastic Navier-Stokes equations drive L\'evy noises. J. Funct. Anal. 272(1), 227-254 (2017).

%\bibitem{Fre} Freidlin, M.I.: Random perturbations of
%reaction-diffusion equations: The quasi-deterministic approach. Trans. Amer. Math. Soc. 305, 665-697 (1988)

%\bibitem{FW} Freidlin, M.I., Wentzell, A.D.: Random perturbation of dynamical systems. Translated by Szuc, J. \ Springer. Berlin (1984)

\bibitem{FS}   Foondun M, Setayeshgar L.:  Large deviations for a class of semilinear stochastic partial differential equations.   Statist. Probab. Lett.  {121},  143-151 (2016).


\bibitem{HS} Hall P., Schimek M.: Moderate-deviation-based inference for random degeneration in paired rank lists,
J. Amer. Statist. Assoc. 107, 661-672(2012).

\bibitem{Gao}  Gao F.: Moderate deviations for a nonparametric estimator of sample coverage, Ann. Statist. 41, 641-669 (2013).

\bibitem{GW} Gao,  F.,  Wang  S.:  Asymptotic behaviors for functionals of random dynamical systems. Stoch. Anal. Appl. {34}(2), 258-277 (2016)

%\bibitem{G} Guillin, A.:  Averaging principle of SDE with small diffusion: moderate deviations.  Ann. Probab. 31, 413-443 (2003)


%\bibitem{GL} A. Guillin, R. Liptser, Examples of moderate deviation principle for diffusion processes.
%Discrete Contin. Dyn. Syst. Ser. B 6 (2006) 803-828


  \bibitem{Gy}Gy\"ongy  I.: Existence and uniqueness results for semilinear stochastic partial differential equations.  Stochastic Process. Appl. 73, 271-299 (1998).
%\bibitem{HS}Hausenblas, E., Seidler, J.: A note on maximal inequality for stochastic convolutions.  Czechoslovak Mathematical Journal,  51(4), 785-790 (2001)
\bibitem{Ich} Ichikawa  A.: Some inequalities for martingales and stochastic convolutions. Stochastic Anal. Appl. 4, 329-339 (1986).


 %\bibitem{KX} Kallianpur, G., Xiong, J.: Large deviations for a class of stochastic partial differential equations. Ann. Probab. 24, 320-345 (1996)
\bibitem{LWYZ}
Li Y., Wang R., Yao N.,  Zhang, S.: A moderate deviation principle for stochastic Volterra equation. Statist.  Probab. Lett., 122(10), 79-85 (2017).


\bibitem{Klebaner1999} Klebaner F., Liptser R.:  Moderate deviations for randomly
perturbed dynamical systems, Stochastic process. Appl. 80, 157-176 (1999).


 %\bibitem{MiaoShen} Y. Miao, S. Shen, Moderate deviation principle for autoregressive processes.  J. Multivariate Anal.  100 (2009) 1952-1961.



\bibitem{Walsh} Walsh  J.: An introduction to stochastic partial differential equations, in P.L. Hennequin (eds.),
 \'{E}cole d'\'{e}t\'{e} de Probabilit\'{e}s St. Flour XIV, in: Lect. Notes Math., vol. 1180, Springer, Berlin, Heidelberg, New York  (1986)

\bibitem{WZZ}  Wang  R., Zhai  J.,  Zhang, T.:  A moderate deviation principle for 2-D stochastic Navier-Stokes equations. J.  Differential Equations {258}, 3363-3390 (2015).

\bibitem{WZ} Wang R., Zhang T.: Moderate deviations for stochastic reaction-diffusion equations with multiplicative noise.  Potential Anal.   {42},  99-113 (2015).


    % \bibitem{Wu} Wu, L.: Moderate deviations of dependent random variables related to CLT. Ann. Probab. 23, 420-445 (1995)

%\bibitem{XZ1} Xu, T., Zhang, T.: Large deviation principles for 2-D stochastic Navier-Stokesequations driven by L\'{e}vy processes. J. Funct. Anal.  257, 1519-1545 (2009)


\bibitem{YGD} Ye H., Gao  J., Ding  Y.: A generalized Gronwall inequality and its application to a fractional differential equation. J. Math. Anal. Appl. 328, 1075-1081 (2007).

\end{thebibliography}
\end{document}